\documentclass[12pt]{amsart}
\usepackage{amsmath,amsthm,amsfonts,amssymb}
\begin{document} 
\newcommand{\B}{{\mathbb B}}
\newcommand{\C}{{\mathbb C}}
\newcommand{\N}{{\mathbb N}}
\newcommand{\Q}{{\mathbb Q}}
\newcommand{\Z}{{\mathbb Z}}
\renewcommand{\P}{{\mathbb P}}
\newcommand{\R}{{\mathbb R}}
\newcommand{\Sc}{{S_{\C}}}
\newcommand{\rc}{\subset}
\newcommand{\rank}{\mathop{rank}}
\newcommand{\trace}{\mathop{tr}}
\newcommand{\dimc}{\mathop{dim}_{\C}}
\newcommand{\dimr}{\mathop{dim}_{\R}}
\newcommand{\codimr}{\mathop{codim}_{\R}}
\newcommand{\Lie}{\mathop{Lie}}
\newcommand{\Aut}{\mathop{{\rm Aut}}}
\newcommand{\Auto}{\mathop{{\rm Aut}_{\mathcal O}}}
\newcommand{\alg}[1]{{\mathbf #1}}
\newcommand{\tensor}{\otimes}
\newcommand{\Gc}{G_{\C}}
\newtheorem*{definition}{Definition}
\newtheorem*{claim}{Claim}
\newtheorem{corollary}{Corollary}
\newtheorem*{Conjecture}{Conjecture}
\newtheorem*{SpecAss}{Special Assumptions}
\newtheorem{example}{Example}
\newtheorem*{remark}{Remark}
\newtheorem*{observation}{Observation}
\newtheorem*{fact}{Fact}
\newtheorem*{remarks}{Remarks}
\newtheorem{lemma}{Lemma}
\newtheorem{proposition}{Proposition}
\newtheorem{theorem}{Theorem}
\title[Connected Lie Groups As Automorphism Groups]{%
Realizing Connected Lie Groups As Automorphism Groups
Of Complex Manifolds
}
\author {J\"org Winkelmann}
\begin{abstract}
We show that every connected real Lie group can be realized
as the full automorphism group of a Stein hyperbolic
complex manifold.
\end{abstract}
\address{%
J\"org Winkelmann \\
 Korea Institute for Advanced Study\\
 School of Mathematics\\
 207-43 Cheongryangri-dong\\
 Dongdaemun-gu\\
 Seoul\\
 130-012 Korea
}
\email{jwinkel@member.ams.org\newline\indent{\itshape Webpage: }%
http://www.math.unibas.ch/\~{ }winkel/
}
\maketitle
\section{Introduction}
Saerens and Zame, and independently Bedford and Dadok proved that,
given a compact real Lie group $K$ there always exists a strictly
pseudoconvex bounded domain $D\subset\C^n$ such that $Aut(D)\simeq K$.
By the theorem of Wong-Rosay (which states that
every strictly pseudoconvex bounded domain with
non-compact automorphism group is isomorphic to the ball) it is clear
that an arbitrary non-compact real Lie group can not be realized as
the automorphism of a strictly pseudoconvex bounded domain in $\C^n$. 
However, as we proved in an earlier paper \cite{W1}, 
for any connected real Lie
group $G$ there does exist a complex manifold $X$ on which $G$ acts effectively.
Moreover, $X$ can be chosen in such a way that it enjoys several
of the key properties of strictly pseudoconvex bounded domains.
Namely, $X$ can be chosen such that it is both Stein and hyperbolic
in the sense of Kobayashi.

The purpose of the present note is to prove that it is possible to
rule out additional automorphisms, i.e.~it is possible to achieve
$Aut(X)\simeq G$. 
\begin{theorem}
Let $G$ be a connected real Lie group.
Then there exists a Stein, complete hyperbolic complex manifold 
$X$ on which $G$
acts effectively, freely, properly and with totally real orbits
such that $\Aut(X)\simeq G$.
\end{theorem}
The idea is to follow the strategy of Saerens and
Zame:
Construct the desired manifold as an open subset
of a larger Stein manifold in such a way that the given group
acts on this open subset.
Ensure that every automorphism of this open subset
can be extended to the boundary,
then
modify the boundary in such a way that this $CR$-hypersurface simply
has no automorphisms other than those from the given group. 
The latter can be done using the fact that a
$CR$-hypersurface (unlike a complex manifold) does have local
invariants.
A principal difficulty in this approach is to obtain an extension 
of automorphisms of the open subset to
the boundary. 
If one is concerned only with compact Lie groups,
then one can work with a strictly pseudoconvex
bounded domain $D$.
For such a domain it is evident that 
for every automorphism $\phi$ of $D$
there exists a sequence  
$x_n\in D$ such that both $x_n$ and $\phi(x_n)$ converge to a strictly
pseudoconvex point in
the boundary. 
This is the starting point for the extension of the automorphism $\phi$
to the boundary $\partial D$.

Now, our goal is to obtain a result for arbitrary connected Lie groups,
which are not necessarily compact.

This lack of compactness assumption creates
some difficulties.

There are two main problems:
First, an arbitrary non-compact Lie group is not necessarily
linear. For instance, the universal cover of $SL_2(\R)$ cannot
be embedded into a linear group.
Second, as already mentioned, the theorem of Wong-Rosay implies that
in general a non-compact Lie group can not
be realized as the full automorphism group
of a strictly pseudoconvex bounded domain with
smooth boundary.
Thus we have to work with domains which are not bounded or
where the boundary is not everywhere smooth.
The trouble is that it is therefore no longer clear that
for every automorphism $\phi$ there exists a sequence $x_n$
in the domain such that both $x_n$ and $\phi(x_n)$ converge
to a nice point in the boundary.

In \cite{TS} a result similar to ours
is claimed for certain Lie groups with a rather sketchy outline
of a possible proof.

The first of the aforementioned two problems 
is dealt with by {\sl assuming}
the group $G$ to be linear while the second problem is simply ignored.
Since the second problem is in fact a serious obstacle, the proof
sketched in \cite{TS} can not be regarded as complete.

We proceed in the following way:
To deal with the first problem, we note that every Lie algebra
is linear by the theorem of Ado. Therefore, in a certain sense,
every Lie group is linear up to coverings and the first problem
can be attacked by working carefully with coverings.

For the second problem, we use bounded domains whose boundaries 
are smooth
outside an exceptional set $E$ which is small in a certain sense.
Exploiting this smallness we
prove that for every automorphism $\phi$ there must
exist a sequence $x_n$ such that both $x_n$ and $\phi(x_n)$
converge to a boundary point outside the ``bad set'' $E$.

Once this has been verified, we can prove 
(using arguments similar to those used in\cite{SZ},
\cite{BD}) that $\phi$ extends
as holomorphic map near $\lim(x_n)$, and use the theory
of Chern-Moser-invariants
to deduce that $\phi$ was in fact given by left multiplication
with an element of $G$.

\subsection{Disconnected Lie Groups}
The result of Bedford and Dadok resp.~Saerens and Zame is
valid for all compact groups, not only connected ones.
However, compactness implies that in this case there are no more
than finitely many connected components.

We conjecture that  our main theorem is valid for arbitrary real
Lie groups, including those with finitely or countably infinitely
many connected components.

As a first step regarding disconnected Lie groups,
we proved in \cite{W2} that the statement of our main theorem does
hold for countable discrete groups.

\section{Linearization}
Given a real Lie group $G$, we look for a bounded domain on which
this group acts. For this purpose we use the theory of hermitian
symmetric spaces.

We will need the following:
\begin{proposition}\label{splinear}
Let $\tilde G$ be a simply-connected real Lie group.

Then there exists a natural number $n$ and a  Lie group homomorphism
$\xi:\tilde G \to Sp(2n,\R)$
such that the following conditions are fulfilled:
\begin{enumerate}
\item
$\xi$ has discrete fibers.
\item
The image $\xi(\tilde G)$ is closed in $Sp(2n,\R)$.
\end{enumerate}

\end{proposition}
\begin{proof}
By Ado's theorem there is an injective Lie algebra homomorphism
$\Lie(\tilde G)\to \Lie GL(m,\R)$ for some $m\in\N$.
Since $\tilde G$ is simply-connected, this induces a Lie group
homomorphism $\xi_0:\tilde G\to GL(m,\R)$ with discrete fibers.
Let $V=\R^m$ and $W=V\oplus V^*$ where $V^*$ is the 
vector space dual of $V$. Then $W$ carries a natural
symplectic structure given by 
\[
(v,\lambda)\cdot(v',\lambda')= \lambda(v')-\lambda'(v)
\]
which is evidently preserved by the natural diagonal
action of $GL(V)$ on $W$.
Hence there is an embedding $i:GL(m,\R)\hookrightarrow Sp(2m,\R)$.

Let $\xi_1=i\circ \xi_0:\tilde G\to Sp(2m,\R)$, $H=\xi_1(\tilde G)$ and
$H'$ its commutator group.
Then $H'$ is already closed in $Sp(2m,\R)$. The quotient group
$H/H'$ is a connected
commutative real Lie group, hence $H/H'\simeq(S^1)^k\times(\R)^l$
for some $k,l\in\N\cup\{0\}$. 
It is easy to see that there is a closed embedding $j:H/H'\hookrightarrow
Sp(2m',\R)$ for some $m'\in\N$ .
Furthermore there is an embedding
$\zeta:Sp(2m,\R)\times Sp(2m',\R)\hookrightarrow Sp(2n,\R)$
with $n=m+m'$.
Now let $\tau:H\to H/H'$ denote the natural projection
and define $\xi:\tilde G\to Sp(2n,\R)$ by
\[
\xi(g)=\zeta(\xi_1(g),j(\tau(\xi_1(g)))).
\]
\end{proof}

\section{Hermitian symmetric domains}
For basic facts on symmetric spaces, see e.g.~\cite{H}.

Let $S=Sp(2n,\R)$  and let $K$ denote a maximal compact subgroup.
Then the quotient manifold $D_0=S/K$ can be endowed with the
structure of a hermitian symmetric domain.
Furthermore there exist open embeddings (``Cayley transform'')
\[
D_0\hookrightarrow \C^N \hookrightarrow Q
\]
such that 
\begin{enumerate}
\item
$D_0$ is relatively compact in $\C^N$,
\item
$Q$ is a projective manifold (the ``compact dual of D'') and
\item
the $Sp(2n,\R)$-action on $D_0$ extends to an 
$Sp(2n,\C)$-action on $Q$.
\end{enumerate}

\begin{lemma}\label{free}
Let $Q$ be a complex manifold on which a complex Lie group
$\Sc$ acts holomorphically and $D_0\subset Q$ a non-empty open subset.

Then there exists a natural number $m\in\N$ and points
$p_1,\ldots,p_m\in D_0$ such that
\[
\cap_{i=1}^m \{ g\in \Sc: g(p_i)=p_i \} =\{e\}
\]
\end{lemma}
\begin{proof}
We choose a sequence of points $p_i\in D_0$ recursively.
First $p_1$ is chosen arbitrarily.
When $p_1,\ldots,p_k$ are already chosen,
we define $I_k=\{g\in \Sc:g(p_i)=p_i,\ 1\le i\le k\}$.
Then we proceed as follows:
If $\dim(I_k)>0$, we choose $p_{k+1}$ such that there is an
element $a_{k+1}$ in the connected component $I_k^0$
such that $a_{k+1}(p_{k+1})\ne p_{k+1}$. This ensures
$\dim I_{k+1}<\dim I_k$.
If $\dim I_k=0$, then $I_k$ is countable.
Thus 
\[
\Lambda=\cup_{g\in I_k\setminus\{e\}}\{x\in Q:g(x)=x\}
\]
is a countable union of nowhere dense analytic subsets of $Q$.
It follows that $\Lambda$ is a set of measure zero for any Lebesgue
class measure on $Q$. In particular $\Lambda\cap D_0\ne D_0$ and we can choose
$p_{k+1}\in D_0\setminus\Lambda$. By the definition of $\Lambda$
this choice enforces $I_{k+1}=\{e\}$.
\end{proof}

\begin{proposition}\label{prop3}
Let $\tilde G$ be a simply-connected real Lie group.
Then there exists a discrete central subgroup $\Gamma$
such that for $G=\tilde G/\Gamma$ the following properties hold:

There exists a natural number $N$, a
bounded domain $D\subset\C^N$, complex analytic subsets
$E\subsetneq \C^N$, $Z\subset D$ and a $G$-action on $Z$ such that
\begin{enumerate}
\item
There is a $G$-invariant non-empty open subset $\Omega$ of $Z$
such that $G$ acts freely, properly, and with totally real orbits
on $\Omega$.
\item
The topological closure $\bar \Omega$ of $\Omega$ in $\C^N$ is contained
inside $Z\cup E$.
\item
$\Omega\cap E=\{\}$.
\item
$\dim_{\C}(\Omega)\ge 3$.
\end{enumerate}
\end{proposition}
\begin{proof}
By prop.~\ref{splinear}, 
there is a discrete central subgroup $\Gamma$ of $\tilde G$ 
such that $G=\tilde G/\Gamma$ can be embedded into some $Sp(2n,\R)$
as closed Lie subgroup.
Let $D_0=Sp(2n,\R)/K$ be the associated hermitian symmetric space
and $Q$ and $\Sc=Sp(2n,\C)$ as described in the beginning of this
section. 

By lemma~\ref{free} there is a natural number
$m$ and a point $p\in D=D_0^m$ such that 
the diagonal $\Sc$-orbit in $Q^m$ through $p$ is free.
Now let $\bar E=\{ x\in Q^m: \dim \Sc(x)<\dim \Sc\}$,
$Z'=\Sc(p)$ and $Z=Z'\cap D$.
Because the $\Sc$-action on $Q^m$ is algebraic, 
the $\Sc$-action on $Q^m$ is algebraic as well.
In particular every $\Sc$-orbit
in $Q^m$ is Zariski open in its closure.
This implies in particular that $\overline{Z'}\subset Z'\cup E$.

Now $G$ is closed in $Sp(2n,\R)$ and $Sp(2n,\R)$ is closed in
$Sp(2n,\C)=\Sc$. 
We obtain a fiber bundle
$\tau:\Sc\to G\backslash\Sc$,
where $G\backslash\Sc$ denotes the quotient of $\Sc$
by the left action of $G$.
Let $U\subset G\backslash\Sc$
be a relatively compact open contractible subset
and $\Omega=\{g\cdot p: g \in \tau^{-1}(U)\}$.

Then $\Omega$ has the desired properties.
(Concerning property $(4)$, observe that 
$\dimc(\Omega)=\dimc(\Sc)\ge 3$ by our construction.)
\end{proof}

\section{Chern-Moser-invariants}
\subsection{Chern-Moser-invariants}
For every real-analytic strictly pseudoconvex $CR$-hypersurface $M$
in a complex manifold $X$
and every point $p\in M$ there is a system of local coordinates 
\[
(w;z)=(w;z_1,\ldots,z_n)
\]
($w,z_i\in\C$, $n+1=\dimc(X)$)
such that $M$ can be written as $M=\{\rho<0\}$
where $\rho$ is a real-analytic function
whose power series development is given as
\[
\rho(w;z)=\Im(w) + ||z||^2 +\sum_{k,l\ge 2} F_{k,l,r}(\Re(w),z,\bar z)
\]
where $F_{k,l,r}$ is a polynomial of bidegree $(k,l)$ in $z$ and $\bar z$
and degree $r$ in $\Re(w)$.

A point $p\in M$ is called {\em umbilical} if $F_{2,2,0}=0$.
For non-umbilical points we define scalar invariants $K_{k,l,r}$
(for $k,l\ge 2, r\in\N$)
given by $K_{k,l,r}=||F_{k,l,r}||^2$ where $||\ ||$ denotes the
euclidean norm, i.e., the norm induced by the scalar
product for which the monomials in the coordinates
$\Re(w),z_i,\bar z_i$ constitute an orthonormal basis.

If $x,y$ are non-umbilical points on $M$ such that the $CR$-hypersurface
germs $(M,x)$ and $(M,y)$ are isomorphic, then
all these invariant $K_{k,l,r}$ must assume the same values at $x$ and $y$. 

For convenient application later on, we define
$K_d=\sum_{k+l=d}K_{k,l,0}$ for $d\ge 4$.
\subsection{Jet bundles}
We recall the notion of jets (see \cite{GG}):
For manifolds $X$ and $Y$ and points $x\in X$, $y\in Y$,
the set of $k$-jets $J^{k}(X,Y)_{x,y}$ is the set of
equivalence classes of 
map germs where two real-analytic 
map germs are equivalent iff their respective
Taylor series developments agrees up to order $k$.
$J^k(X,Y)$ is the disjoint union of all $J^{k}(X,Y)_{x,y}$
(with $x\in X$ and $y\in Y$). There is a natural manifold structure
on $J^k(X,Y)$ for which we obtain a fiber bundle 
({\em``source map''}) $\alpha:J^k(X,Y)\to X$.

\subsection{Transversality}
We will need the {\em multijet transversality theorem}
(\cite{GG}, thm.~4.13).
Let $X^{(s)}$ denote the space of those $s$-tuples $(x_1,\ldots,x_s)\in X^s$
where the $x_i$ are all {\em distinct} elements in $X$.
Let
\begin{multline*}
J^k_s(X,Y)=\{(f_1,\ldots,f_s)\in \left(J^k(X,Y)\right)^s:
\alpha(f_1),\ldots,\alpha(f_s)\in X^{(s)}\}
\end{multline*}

Then each $f\in C^\infty(X,Y)$ induces a map
$j^k_s(f):X^{(s)}\to J^k_s(X,Y)$ in a natural way.

Let $W$ be a submanifold of 
codimension $c$ in $J^k_s(X,Y)$.

Then the {\em multijet transversality theorem} implies that
the function space $C^\infty(X,Y)$ contains a {\em residual} subset
$A$ such that $\left(j^k_s(f)\right)^{-1}(W)$ is of codimension at least
$c$ in $X^{(s)}$.

\begin{remark}
\begin{enumerate}
\item
In the statement on the codimension,
the codimension of the empty set is to be understood as $+\infty$.
\item
A subset of a topological space $V$ is called {\em residual}
if it is the intersection of countably many open dense subsets.
If $V$ has the {\em Baire property}, then every residual subset
of $V$ is dense. The function spaces $C^\infty(X,Y)$ and 
$C^\omega(X,Y)$ have the Baire property (for any pair of
manifolds $(X,Y)$.)
\item
Similar results hold for the function spaces of type $C^{\omega}$,
i.e.~real-analytic mappings, which in fact can be deduced from
the transversality results for $C^{\infty}$-maps,
using the fact that $C^\omega$-maps are dense in $C^\infty$.
\item 
In the real-analytic category, $W$ does not need to be smooth, it suffices
if $W$ is a (possibly singular) real-analytic subset.
As explained in \cite{SZ}, this can be verified using the fact that a
real analytic subset $W$ admits a stratification 
$W=W_0\supset W_1 \supset W_2\ldots$ such that each $W_k\setminus W_{k+1}$
is smooth.
\end{enumerate}
\end{remark}

\subsection{A proposition}

Let us now assume that there is a real Lie group $G$ acting holomorphically
on $X$ with totally real orbits. 
Let us furthermore assume that the action in proper.
Then orbits can be separated by invariant functions.
Around any given point $p\in X$, we may choose
local holomorphic coordinates $x_i$
in such a way that $x_i(p)=0\ \forall i$ and 
\[
T_p\left(G\cdot p\right)\subset T_p\left(\{x:\Re( x_i)=0\ \forall i\}\right)
\]
It follows that for a every real homogeneous polynomial $P$
of degree $k$ there is a $G$-invariant real-analytic function $f$ 
defined on some open neighbourhood of $p$ in $X$
such that $P(x_1,\ldots,x_n)=f(x)+O(||x||^{k+1})$.
As a consequence, we obtain the statement below:
\begin{lemma}\label{lemma-X}
Let $G$ be a real Lie group acting holomorphically
and properly on a complex manifold $X$ with totally real orbits. 
Assume $\dim_{\C}(X)\ge 2$. 

For $x\in X$, $t\in\R$ let $J^k_+(X,\R)^G_{x,t}$ denote the set of
all $k$-jets induced by germs of $G$-invariant functions $f$
for which the $CR$-hypersurface germ defined by $f=t$ is strictly
pseudoconvex around $x$.
Let $J^k_+(X,\R)^G=\cup_{x\in X,t\in\R}J^k_+(X,\R)^G_{x,t}$.

Then $K_4=\ldots =K_k=0$ defines a real-analytic subspace of codimension
at least $k-3$ in $J^k_+(X,\R)^G$.
\end{lemma}

Now we can prove the proposition given below.

\begin{proposition}\label{CM}
Let $G$  be a real Lie group acting holomorphically
and properly on a complex manifold $X$ with totally real orbits.
Assume $\dim_{\C}(X)\ge 2$. 
Let $p\in X$.

Then $G\cdot p$ admits an open $G$-invariant neighbourhood
$\Omega$ such that:
\begin{enumerate}
\item
The inclusion map $G\cdot p \hookrightarrow \Omega$ is a homotopy
equivalence.
\item
The boundary $\partial \Omega$ is everywhere smooth, real-analytic
and strictly pseudoconvex.
\item
There exists a nowhere dense real-analytic subset $\Sigma\subset
\partial\Omega$ such that for every $x,y\in\partial\Omega\setminus\Sigma$
the $CR$-hypersurface germs $(\partial\Omega,x)$, $(\partial\Omega,y)$
are isomorphic if and only if $x=g\cdot y
$ for some $g\in G$.
\end{enumerate}
\end{proposition}

\begin{proof}
Let $r=\dimr(X)-\dimr(G)$, $B$ the open unit ball in $\R^r$ and
$i:B\hookrightarrow X$ a real-analytic embedding 
with $i(0)=p$
which is everywhere transversal to the
$G$-orbits.
Then $W=G\cdot i(B)$ is an open $G$-invariant neighbourhood
of the $G$-orbit $G\cdot p$.
Since the $G$-action on $X$ is free and proper, we may and do assume
that the map $G\times B\to W$ given by $(g,x)\mapsto g\cdot i(x)$
is bijective.

We define $\rho_0\in C^\omega(W)$ via $\rho_0(g\cdot i(v))
=||v||^2$ for $g\in G$, $v\in B$.

An easy calculation in local coordinates shows that
\[
W_\epsilon=\{x\in W:\rho_0(x)<\epsilon \}
\]
is strictly pseudoconvex for all sufficiently small $\epsilon>0$.
We fix now a number $1>\delta>0$ such that $W_\delta$
is strictly pseudoconvex.

Then $W_\delta$ is a $G$-invariant
open neighbourhood of $G\cdot p$ fulfilling
conditions $(1)$ and $(2)$  of the proposition.
To achieve condition $(3)$, we have to modify the
defining function $\rho_0-\delta$ of the  open domain $W_\delta$
using the theory of Chern-Moser invariants.

Every function on $B\simeq i(B)$ extends uniquely to a
$G$-invariant function on $W$; this yields a 
bijective map 
\[
\zeta: C^\omega(B)\longrightarrow C^\omega(W)^G.
\]

Now let $\Theta$ be an open neighborhood of $(\zeta^{-1}(\rho_0),\delta)$
in $ C^{\omega}(B)\times\R$ such that the following properties
hold for all 
$(f,t)\in\Theta$:

\begin{enumerate}
\item
$\{v\in B: f(v)<t\}$ is a contractible relatively compact open
subset with smooth boundary in $B$;
\item
The domain 
\[
\{x\in W:\zeta(f-t)(x)<0\}=\{g\cdot i(v):f(v)<t,v\in B,g\in G\}
\]
is everywhere strictly pseudoconvex.
\end{enumerate}

Let $J^k_+(B,\R)=\cup_{v\in B,t\in\R}J^k_+(B,\R)_{v,t}$ 
where $J^k_+(B,\R)_{v,t}$  denotes the set of
all $k$-jets induced by germs of functions $f:(B,v)\to(\R,t)$
for which the $CR$-hypersurface germ 
\[
\{g\cdot i(x): g\in G, x\in B, f(x)=t\}
\]
is strictly
pseudoconvex around $i(v)$.
For $k\in\N$, $4\le d\le k$ we define functions $\tilde K_d$ on
$J^k_+(B,\R)$ as follows:
If $j$ is the $k$-jet at $v\in B$ for some map germ $f:(B,v)\to (\R,t)$,
then $\tilde K_d(j)$ is defined as the scalar invariant $K_d$
for the $CR$-hypersurface $\{y\in W:(\zeta(f))(y)=t\}$ at $i(v)$.

We define the ``umbilical locus'':
\[
U_k=\{j\in J^k_+(B,\R):\tilde K_4(j)=0\}
\]
and the ``locus of coinciding scalar invariants'':
\[
E_k=\{(j_1,j_2)\in J^k_+(B,\R)^2:\tilde K_d(j_1)=\tilde K_d(j_2)\ 
\forall 4\le d \le k \}
\]

Since $J^d_+(B,\R)$ is an open subset 
in $J^d(B,\R)$, $U_k$ and $E_k$ can be
regarded as locally closed real-analytic subset in $J^d(B,\R)$
resp.~$J_2^d(B,\R)$.

Fix $k$ such that $k-3>2\dimr(B)$.
Then lemma~\ref{lemma-X} implies that the codimension of $E_k$
exceeds the dimension of $B\times B$.

The {\em multijet transversality theorem}
implies that there is a residual set $A\subset C^{\omega}(B,\R)$
such that every $f\in A$ is transversal to both $U_k$ and $E_k$.

Since $A$ is {\em residual}, it is dense
in $C^\omega(B,\R)$.
Therefore $A\times\R$ intersects
the open set $\Theta$. Let $(\rho_1,t_0)\in (A\times\R)\cap\Theta$.
Let $\Sigma_0\subset W$ be the set of all points $x\in W$ such that
the $CR$-hypersurface $\{y\in W:\zeta(\rho_1)(y)=\zeta(\rho_1)(x)\}$ is
umbilical at $x$. Then transversality of $\rho_1$ with 
respect to $U_k$ implies
that $\Sigma_0$ is a nowhere dense, locally closed
real-analytic subset of $W$.
As a consequence, we can find a real number $t$ close to $t_0$
such that $(\rho_1,t)\in\Theta$ and such that
$\Omega=\{y\in W:\zeta(\rho_1)(y)<t\}$ has the following property:

{\em
`` $\Sigma_0\cap\Omega$ is nowhere dense in $\Omega$.''
}

Now $(\rho_1,t)\in\Theta$ implies that conditions $(1)$ and $(2)$
are fulfilled for our choise of $\Omega$.
Furthermore transversality of $\rho_1$ with respect to $E_k$ 
(in combination with $\codimr(E_k)>\dimr(B)$)
implies that $\Omega$ fulfills
condition $(3)$ of the proposition.
This completes the proof.
\end{proof}

\section{Privalov's theorem}
We are now in position to use the classical theorem of Privalov
in order
to show that for every automorphism $\phi$ there is a
sequence $x_n$ such that both $x_n$ and $\phi(x_n)$ converge
to a point in the good part of the boundary.

\begin{proposition}\label{priv}
Let $D$ be a bounded domain in $\C^N$, $E\subset \C^N$, $Z\subset D$
closed analytic subsets, $\Omega$ an open subset of $Z$,
$M$ its boundary in $Z$. Assume $E\cap\Omega=\{\}$.
Assume that $M$ is everywhere smooth and that the closure of $\Omega$
in $\C^N$ is contained in $\Omega\cup M \cup E$.
Let $\Omega'$ be the closure of $\Omega$ in $Z$, i.e., 
$\Omega'=\Omega\cup M$.

Furthermore let $\tilde\Omega$ denote the universal covering of
$\Omega$ and $\pi:\tilde\Omega'\to\Omega'$ and $\tilde M\to M$ the
corresponding coverings.

Then for every holomorphic automorphism $\phi\in\Aut(\tilde\Omega)$
there is a sequence $x_n$ in $\tilde\Omega$ and points 
$q,\bar q\in \tilde M$
such that $\lim x_n=q$ and $\lim\phi(x_n)=\bar q$.
\end{proposition}
\begin{proof}
Fix $\phi\in\Aut(\tilde\Omega)$.
Let $\Delta$ be the unit disk in $\C$, $\bar\Delta$ its closure
in $\C$ and
$\partial\Delta$ its boundary.

We choose a $C^\infty$ map $\zeta:\bar\Delta\to\tilde\Omega'$
such that
\begin{enumerate}
\item
$\zeta|_{\Delta}$ maps $\Delta$ holomorphically
into $\tilde\Omega$.
\item
$\zeta^{-1}(\tilde M)$ is a subset of positive Lebesgue measure
in $\partial\Delta\simeq S^1$.
\end{enumerate}
Now we consider $\eta:\Delta\to\C^N$ given by $\eta=\pi\circ\phi\circ\zeta$.
Then $\eta$ is a $N$-tuple of bounded holomorphic functions.
It follows (\cite{K},\cite{R}) that the non-tangential limit
exists almost everywhere on $\partial\Delta$.
For $t\in\partial\Delta$, let $\lim_{n-t}\eta(t)$ denote
this non-tangential limit.
Evidently $\lim_{n-t}\eta(t)\in\Omega'\cup E$ whereever defined.
We claim that $A=\{t:\lim_{n-t}\eta(t)\in E\}$
is a set of measure zero.
Indeed $t\in A$ implies that for every holomorphic function
$f$ on $\C^N$ which vanishes on $E$, we obtain
\[
\lim_{n-t}(f\circ\eta)(t)=0
\]
If $A$ is not a set of measure zero, it would follow
from Privalov's theorem (\cite{K}) that
$f\circ\eta$ would vanish for every such $f$.
But this would imply $\eta(\Delta)\subset E$,
contradicting $\eta(\Delta)\subset \Omega$.
Thus $A$ must be a set of measure zero.
It follows that there exists a point $q\in\partial\Delta\cap
\zeta^{-1}(\tilde M)$
such that the non-tangential limit for $\eta$ exists at $q$
and is not in $E$.

Now fix a triangle $T\subset\Bar\Delta$ with its three edges on $\partial D$ 
one of which is $q$
($T$ denotes the triangle with interior, i.e., the convex hull
spanned by the three edges.)
By the definition of the notion ``non-tangential limit'' we
have a limit
\[
\lim_{x\in T,x\to q}\eta(x)=v\in \Omega'\subset\C^N
\]
and thus a continuous map $\bar\eta:T\cup\{q\}\to\Omega'$ with
$\bar\eta|_T=\eta$.  Let $W$ be a simply-connected open neighbourhood of
$v$ in $\Omega'$, and $V$ an open connected neighbourhood of $q$ in
$\bar\eta^{-1}(W)$.  Observe that $\pi:\tilde\Omega'\to\Omega'$ is an
unramified covering.  
Since $W$ is simply-connected, it follows that $\pi^{-1}(W)$ is a disjoint
union of connected components each of which is isomorphic to $W$.
Connectedness of $V$ implies that $\phi(\zeta(V))$ is
contained in one connected component of $\pi^{-1}(W)$. Together with
$\lim_{x\in T,x\to q}\eta(x)=v$ this implies that there is a point
$\tilde v\in\pi^{-1}(v)$ such that
\[
\lim_{x\in T,x\to q}\phi(\zeta(x))=\tilde v=\bar q
\]

For any sequence $t_n$ in $int(T)$ converging to $t$
we now obtain a sequence $x_n=\zeta(t_n)$ with
convergent limits $\lim x_n=q\in\tilde M$, 
$\lim\phi(x_n)=\bar q\in\tilde\Omega'$.

Finally we note that $\bar q$ cannot be in $\tilde\Omega$: 
$\phi$ is an automorphism of $\tilde\Omega$ and
therefore $\lim x_n\not\in\tilde\Omega$ implies that
$\phi(x_n)$ cannot converge inside of $\tilde\Omega$ . 
Hence $\bar q\in\tilde M$.
\end{proof}

\section{Extension through the boundary}

We need the following well-known extension result.

\begin{proposition}\label{ext}
Let $\Omega$ be an open subset in a Stein manifold $Z$.
Assume that there are points $q,\bar q\in\partial\Omega$, an automorphism
$\phi\in\Aut(\Omega)$, and a sequence of points $x_n\in\Omega$
with $\lim x_n=q$ and $\lim\phi(x_n)=\bar q$.
Assume in addition that $\partial\Omega$ is real-analytic and
strictly pseudoconvex near $q$ and $\bar q$.

Then there exists an open neighbourhood $V$ of $q$ in $Z$
and a holomorphic map $\Phi:V\to Z$ such that 
$\Phi|_{\Omega\cap V}=\phi|_{\Omega\cap V}$.
\end{proposition}
\begin{proof}
First, \cite{FR} implies that $\phi$ can be extended to a continuous
map $\bar\phi$ on $\bar\Omega$ near $q$. 
Since $\bar\phi$ is continuous and $\bar\phi|_\Omega$ is holomorphic,
it is clear that $\bar\phi|_{\partial\Omega}$ is a continuous $CR$-map.
(For a not necessarily differentiable function the notion ``$CR$-map''
is defined via regarding derivatives in the sense of distributions.
Then the condition ``$CR$'' translations into the vanishing of certain
integrals involving test functions -- a closed condition; hence
holomorphy of $\bar\phi|_{\Omega}=\phi$ implies that
$\bar\phi|_{\partial\Omega}$ is a $CR$-map.)

Thus \cite{Be} implies that this extension
is already $C^\infty$ and finally \cite{BR} or \cite{DF}
yield that there is a holomorphic extension into some open neighbourhood.
\end{proof}

\section{Rigidity}
\begin{lemma}\label{real-valued}
Let $\Omega$ be a strictly pseudoconvex domain in a Stein manifold $V$.
Let $f$ be a holomorphic function on $V$ such that $f(\partial\Omega)
\subset\R$.

Then $f$ is constant.
\end{lemma}
\begin{proof}
By the assumption of $\Omega$ being strictly pseudoconvex
it follows that for every point $p\in\Omega$ close enough
to the boundary there exists a continuous map $\zeta:\bar\Delta\to V$
such that
\begin{enumerate}
\item $\zeta$ is holomorphic on $\Delta$,
\item $\zeta(0)=p$,
\item $\zeta(\partial\Delta)\subset\partial\Omega$
\end{enumerate}
Now the maximum principle applied to the plurisubharmonic
function $g(x)=\left(\Im f(x)\right)^2$ implies that $\Im f(p)=0$.
Thus the real-analytic function $\Im f$ vanishes in some open subset
of $\Omega$ and therefore (by identity principle) it vanishes everywhere.
Hence $f$ is both holomorphic and everywhere real-valued and
therefore constant.
\end{proof}

\begin{proposition}\label{rigid}
Let $\Omega$ be an open $G$-invariant subset of a complex
manifold $Z$ on which $G$ acts freely with totally real
orbits.
Assume that the boundary $\partial\Omega$ is a smooth $CR$-hypersurface.

Let $\phi$ be an automorphism of $\Omega$, $q\in\partial\Omega$
and $V$ an open neighbourhood of $q$ in $Z$ such that
$\phi|_{V\cap\Omega}$ extends to a holomorphic map $\bar\phi:V\to Z$.

Assume that for every $x\in V\cap\partial\Omega$ both $x$ and 
$\bar\phi(x)$ are contained in the same $G$-orbit.

Assume furthermore that $\partial\Omega$ is strictly pseudoconvex
near $q$.

Then there exists an element $g\in G$ such that $g\cdot x=\phi(x)$
for all $x\in\Omega$.
\end{proposition}
\begin{proof}
Let $g_0\in G$ be such that $\bar\phi(q)=g_0\cdot q$.
We may now replace $\phi$ by the automorphism 
$x\mapsto g_0^{-1}\cdot\phi(x)$
and thereby assume that $\bar\phi(q)=q$.
Now we have to show that $\phi=id_\Omega$.

Let $n=\dim_{\C}(\Omega)$ and $d=\dim_{\R}(G)$.
Let $i:B_{n-d}=\{v\in\C^{n-d}:||v||<1\}\to Z$
be an embedding such that $i(0)=q$ and that $i(B_{n-d})$
is everywhere transversal to the $G$-orbits.
The $G$-action induces a real-analytic map 
$\psi:\Lie(G)\times Z\to Z$
given by $\psi(v,x)=\exp(v)\cdot x$.
This extends to a holomorphic map $\psi_{\C}:U\to Z$
where $U$ is an open neighbourhood of $(0,q)$ in
$(\Lie(G)\tensor\C)\times Z$.
By appropriately shrinking $V$ and $U$ we may assume that $U=N\times V$
where $N$ is an open neighbourhood of $0$ in $\Lie(G)\tensor\C$.
Now we obtain a holomorphic map $\zeta:B_{n-d}\times N\to Z$
via $\zeta(w,v)=\psi_{\C}(v,i(w))$.
Since $B_{n-d}\times N$ is an open domain in
\[
\C^{n-d} \times \Lie(G)\tensor\C \simeq \C^{n-d}\times \C^d \simeq\C^n
\]
this map $\zeta$ yields  local holomorphic coordinates
near $q$.
In these local coordinates
\[
x=(x_1,\ldots,x_n) \mapsto \bar\phi(x)-x
\]
is a holomorphic map all of whose components are real-valued
on $V\cap\partial\Omega$.
Because $\partial\Omega$ is strictly pseudoconvex near $q$,
it follows that this map is constant (lemma~\ref{real-valued}).
Since $\bar\phi(q)=q$, constancy means that it is constant zero.
Thus $\bar\phi\equiv id_V$. Finally, by identity principle it follows
that $\phi(x)=x$ for all $x\in\Omega$, as desired.
\end{proof}

\section{Reduction to the simply-connected case}
\begin{lemma}\label{sc}
Let $G$ be a connected real Lie group, $\tilde G$ its universal
covering and $\Gamma=\pi^{-1}(\{e\})$ where
$\pi:\tilde G\to G$ is the natural projection map.

Assume that there exists a simply-connected
complex manifold $\tilde X$ with $\Aut(\tilde X)
\simeq\tilde G$ such that the $\Gamma$-action on $\tilde X$ is
free and properly discontinuous.

Let  $X=\tilde X/\Gamma$.
Then $\Aut(X)\simeq G$,
\end{lemma}
\begin{proof}
Every automorphism of $X$ lifts to an automorphism of $\tilde X$,
because $\tilde X$ is a universal covering space for $X$.
Therefore the automorphism group of $X$ is isomorphic to $N/\Gamma$
where $N$ denotes the group of all elements of $\Aut(\tilde X)$
which normalize $\Gamma$. But $\Gamma$ is the kernel of a group
homomorphism, hence normal. Thus $N=\tilde G$ and consequently
$\Aut(X)\simeq N/\Gamma =\tilde G/\Gamma \simeq G$.
\end{proof}

It remains to be shown that $\tilde X$ can be constructed in such a 
way that $X=\tilde X/\Gamma$ will be
Stein and completely hyperbolic.
Complete hyperbolicity is easy, since $\tilde X$ being completely
hyperbolic implies that $X$ is completely hyperbolic, too.

The Stein property is more involved, since for an arbitrary
unramified covering $\tilde X\to X$, Steinness of $\tilde X$
does not imply that $X$ is Stein, too.

\begin{proposition}\label{Stein}
Let $G$ be a real Lie group, $\pi:\tilde G\to G$ its universal covering,
$\Gamma=\pi^{-1}(e)$, and
$\tilde X$ a complex manifold on which $\tilde G$
acts properly and freely with totally real orbits.

Let $p\in\tilde X$.

Then there exists an open $\tilde G$-invariant neighbourhood $\tilde U$ of 
$\tilde G\cdot p$ in $\tilde X$ such that for every
$\tilde G$-invariant locally Stein open  submanifold $\tilde\Omega\subset \tilde U$
the complex quotient manifold $\Omega=\tilde\Omega/\Gamma$ is
Stein.
\end{proposition}
(As usual, $\tilde\Omega\subset \tilde U$ is called locally Stein iff
every point $x\in \tilde U$ admits an open neighbourhood $V$
in $U$ such that $V\cap\tilde\Omega$ is Stein.)
\begin{proof}
Essentially, we follow the argumentation in \cite{W1}.

Let $Z$ denote the center of $\tilde G$. Then there exists
a discrete cocompact subgroup $\Lambda$ 
in $Z$ such that $\Gamma\subset \Lambda$
(\cite{W1}, lemma~1).
Let $G_1=\tilde G/\Lambda$ and $X_1=\tilde X/\Lambda$.

Let $\Gc$ be the simply-connected complex Lie group
corresponding to the complex Lie algebra
$\Lie(G)\tensor\C$ and $j:\tilde G\to\Gc$ the natural Lie
group homomorphism induced by the Lie algebra embedding
$\Lie(G)\hookrightarrow\Lie(G)\tensor\C$.

Let $\psi_0:\Lie(\tilde G)\times\tilde X\to\tilde X$ be the map
induced by the group action via $\psi_0(v,x)=\exp(v)\cdot x$.
Then $\psi_0$ extends to a holomorphic map $\psi$ defined
on some open neighbourhood of $\Lie\tilde G\times\tilde X$
in $\Lie\Gc\times\tilde X$.
This open neighbourhood can be chosen as product $N\times \tilde W$
where $N$ is an open neighbourhood of $\Lie\tilde G$ in $\Lie\Gc$ and 
$W$ is an
open neighbourhood of $p$ in $\tilde X$.

Let $n=\dim_{\C}(X)$ and $d=\dim_{\R}(G)=\dim_{\C}(\Gc)$.
Let $i:B_{n-d}=\{v\in\C^{n-d}:||v||<1\}\to W$
be a holomorphic embedding such that $i(0)=p$ and that $i(B_{n-d})$
is everywhere transversal to the $G$-orbits.

We choose a small open neighbourhood $N_1\subset N$ of $0$ in $\Lie\Gc$
such that the map $\zeta:\tilde G\times N_1\times B\to\tilde X$
given by $\zeta:(g,n,x)\mapsto g\cdot\psi(n,z)$
has the property that $\zeta(g,n,z)=\zeta(g',n',z')$ only if
there is an element $v\in\Lie\tilde G$ such that $g'=g\cdot\exp(v)$
and $\exp(-v)\exp(n)=\exp(n')$.
This is possible, because $\tilde G$ acts freely with totally real
orbits.

For $x\in\zeta(\tilde G\times N_1\times B)$
we define $\xi(x)\in Ad(\Gc)$ by $Ad(g\cdot\exp(n))$ if $x=\zeta(g,n,z)$.
Then $\xi$ is a well-defined, holomorphic and $\tilde G$-equivariant
map from an $\tilde G$-invariant open neighbourhood $W_0$ of $p$
to $Ad(\Gc)$. Moreover $\xi$ is constant along the orbits of the
center $Z$ of $\tilde G$.
Therefore it induces a holomorphic map $\xi_1:W_1\to Ad(\Gc)$
where $W_1$ is the image of $W_0$ under the projection $\tilde X\to X_1$.

Observe that $Ad(\Gc)$ is a linear complex Lie group. It follows
that $Ad(\Gc)$ is Stein (\cite{MM})
and hence admits a strictly plurisubharmonic
exhaustion function $\rho_1:Ad(\Gc)\to\R^+$.

Next we consider the real quotient map 
$\tau:\tilde X\to \tilde X/\tilde G=Y$.
Let $y_1,\ldots,y_r$ be real-analytic  local coordinates on $Y$ with
$y_i(\tau(p))=0$. Then $x\mapsto \sum_i y_i(\tau(x))^2$ defines
a $\tilde G$-invariant
real-analytic function $\rho_0$ on a neighbourhood of $\tilde G\cdot p$,
which is easily verified to be strictly plurisubharmonic 
near $\tilde G\cdot p$.

By appropriately shrinking $W_0$ and $W_1$
we may assume that 
there is an $\epsilon>0$ such that $W_0=\{x:\rho_0(x)<\epsilon\}$.

We reparametrize this function via
\[
\rho_0'(x)=\tan\left(\frac{\pi}{2\epsilon}\rho_0(\tau(x))\right)
\]
Now $\rho_0'\to+\infty$ whenever $\rho_0\to\epsilon$.
Thus $\rho_0'$ is an ``exhaustion function modulo $\tilde G$'',
i.e. it is a $\tilde G$-invariant function which induces a
proper continuous map from $W_0/\tilde G$ to $\R^+$.

Moreover, $\rho_0'$ is strictly plurisubharmonic, because $\tan$
is convex and $\rho_0$ is strictly plurisubharmonic.

Next we recall that by lemma~2 in \cite{W1} the natural map
$W_1\to Y_1\simeq W_1/G \times Ad(\Gc)$ is proper.

Therefore $\rho_1+\rho_0'$ is a continuous exhaustion function on $W_1$.
On the other hand, this function is also strictly plurisubharmonic.
Thus $W_1$ is Stein.

Let $\tilde U=W_0$ and $U=\tilde U/\Gamma$.
Then $U$ is Stein, because $W_1$ is Stein and we have an unramified
covering $U\to W_1$.

Assume that $\Omega$ is a $\tilde G$-invariant open locally Stein submanifold
of $\tilde U$. Then $\Omega/\Gamma$ is a $G$-invariant open submanifold
of $U=\tilde U/\Gamma$ which is evidently locally Stein. 
But locally Stein open submanifolds
of Stein manifolds are Stein. 
Hence $\Omega$ is Stein.
\end{proof}

\section{Proof of the Main theorem}

Here we prove our main theorem.
\begin{proof}
Let $\tilde G$ denote the universal covering of $G$,
$\pi:\tilde G\to G$ the natural projection and $\Gamma=\pi^{-1}\{e\}$.
By prop.~\ref{prop3}
there is a quotient $G_1$ of $\tilde G$ by a 
central discrete subgroup $\Gamma_1$ and
a $G_{1}$-action on a complex manifold $\Omega_1$
which is free, proper and with totally real orbits.
Moreover, there is a number $N$, a bounded domain $D\subset\C^N$ and
closed complex analytic subsets $Z\subset D$, $E\subset \C^N$
and an embedding of $\Omega_1$ as open submanifold in $Z$ such that
the closure of $\Omega_1$ in $\C^N$ is contained in $Z\cup E$ and
$\Omega_1\cap E=\{\}$.

Fix $p\in\Omega_1$. We may replace $\Omega_1$ by some appropriately
choosen invariant open neighbourhood of $G_1\cdot p$.
Therefore we may and do from now on assume  that 
$\pi_1(\Omega_1)\simeq\pi_1(G_1)\simeq\Gamma_1$.
(And we keep this assumption throughout all further replacements
of $\Omega_1$ by invariant open subsets of itself.)
Let $\tilde\Omega$ denote the universal covering of $\Omega_1$.

From prop.~\ref{Stein} we deduce that, after replacing $\Omega_1$
with some $G_{1}$-invariant open subset, we may assume that
$U/\Gamma$ is Stein for every open $\tilde G$-invariant locally Stein
submanifold $U$ of $\tilde\Omega$.

Next we apply prop.~\ref{CM}, again replacing $\Omega_1$ by an appropriate
smaller $G_{1}$-invariant open submanifold. Now $\Omega_1$ has a smooth,
real-analytic and strictly pseudoconvex boundary $B$ in $Z$, and
there is a nowhere dense real-analytic subset $\Sigma\subset B$ such that
for every $x,y\in B\setminus\Sigma$
the $CR$-hypersurface germs $(B,x)$, $(B,y)$
are isomorphic if and only if $x=g\cdot y
$ for some $g\in G_{1}$.

Let $\phi\in Aut(\tilde\Omega)$. 
Let $\tau:\tilde\Omega\to\Omega_1$ be the covering map.
We may assume that there is a $G_1$-invariant open subset $X\subset Z$
such that $\bar\Omega_1\cap Z\subset X$. Furthermore we may assume
that the inclusion $\Omega_1\hookrightarrow X$ induces an isomorphism
of the fundamental groups. Then $\tau:\tilde\Omega\to\Omega_1$ extends
to a covering $\tau':\tilde X\to X$ with $\tilde\Omega\hookrightarrow
\tilde X$.
Let
$\tilde\Sigma$ and $\tilde B$
denote the preimages of $\Sigma$ resp.~$B$ under $\tau'$.
By prop.~\ref{priv} there is a sequence of points $x_n\in\tilde\Omega$
and points $q,\bar q\in\tilde B$ 
such that $\lim x_n=q$ and $\lim\phi(x_n)=\bar q$.
By prop.~\ref{ext} it follows that $\phi$ extends to a holomorphic
map $\Phi$ in an open neighbourhood $U$ of $q$ in $\tilde X$. 
Because $\phi^{-1}$ extends to a holomorphic map near $\bar q$
by the same arguments and $\phi\circ\phi^{-1}=id$,
this extension $\Phi$ is locally biholomorphic 
and $\Phi(\tilde B\cap U)\subset\tilde B$.

Recall that $\tilde\Sigma$ is nowhere dense in $\tilde B$.
Hence there is an element 
$q'\in (\tilde B\cap U)\setminus\left(\tilde\Sigma\cup\Phi^{-1}(\tilde\Sigma)
\right)$.
Upon replacing $q$ by $q'$ and $U$ by 
$U\setminus\left(\tilde\Sigma\cup\Phi^{-1}(\tilde\Sigma)
\right)$,
we may from now on assume that $U\cap\Sigma$ and $\Phi(U)\cap\Sigma$
are both empty.
 
For every $z\in\tilde B\cap U$ the $CR$ hypersurface germs
$(\tilde B,z)$ and $(\tilde B,\Phi(z))$ are isomorphic and consequently
there is an element $g_z\in\tilde G$ such that $g_z\cdot z=\Phi(z)$.

By prop.~\ref{rigid} it follows that there is one element $g\in\tilde G$
such that $\phi(z)=g\cdot z$ for all $z\in\tilde\Omega$.
Thus $\phi\in\tilde G$. Since $\phi$ was an arbitrary automorphism
of $\tilde\Omega$, it follows that $Aut(\tilde\Omega)=\tilde G$.
By lemma~\ref{sc} this implies that $Aut(\Omega)=G$ where 
$\Omega=\tilde\Omega/\Gamma$.

Finally let us discuss the Stein condition and hyperbolicity.
Since $\Omega_1$ injects into a bounded domain $D\subset\C^N$,
it is hyperbolic. Because $\tilde\Omega\to\Omega_1$ and
$\tilde\Omega\to\Omega$ are both unramified coverings,
this implies the hyperbolicity of $\Omega$.
Moreover, by the same arguments as in \cite{W1}, we may conclude
that $\Omega$ is even {\em complete} hyperbolic. 

Concerning the Stein property, let us recall application of
prop.~\ref{Stein} further above.
Our choice of $\Omega_1$ at that time had the property that
$U/\Gamma$ was Stein for every locally Stein open subset $U$ of 
$\tilde\Omega$.
Subsequently we shrank $\Omega_1$, replacing it by some open subset with
strictly pseudoconvex boundary. Clearly an open subset with strictly
pseudoconvex boundary is locally Stein. Therefore $\Omega=\tilde\Omega/\Gamma$
is Stein for our final choice of $\Omega_1$.
\end{proof}


\begin{thebibliography}{Bla}

\bibitem{BR} Baouendi, M.S.; Rothschild, L.P.:
Germs of $CR$-maps between real-analytic hypersurfaces.
\sl Invent. Math. \bf 93\rm, 481-500 (1988)

\bibitem{BD} Bedford, E.; Dadok, J.:
Bounded domains with prescribed group of automorphisms.
\sl Comm. Math. Helv. \bf 62 \rm, 561--572 (1987)

\bibitem{Be} Bell, S.:
Local regularity of CR-homeomorphisms.
\sl Duke Math. J. \bf 57\rm, 295--300 (1988)

\bibitem{Bo} Borel, A.:
Semisimple groups and Riemannian symmetric spaces.
Texts and Readings in Mathematics \bf 16\rm, Hindustan Book Agency.
(1998)

\bibitem{DF} Diederich, K.; Fornaess, J.:
Proper Holomorphic Mappings between real-analytic
pseudoconvex domains in $\C^n$.
\sl Math. Ann. \bf 282\rm, 681-700 (1988)

\bibitem{FR} Forstneric, F.; Rosay, J.P.:
Localization of the Kobayashi Metric and the boundary Countinuity of Proper
Holomorphic Mappings.
\sl Math. Ann. \bf 279 \rm, 239--252 (1987)

\bibitem{F} Frankel, S.:
Complex geometry of convex domains that cover varieties.
\sl Acta Math. \bf 163\rm, 109--149 (1989)

\bibitem{GG} Golubitsky, M.; Guillemin, V.:
Stable Mappings and Their Singularities.
GTM Springer

\bibitem{H} Helgason, S.: Differential geometry, Lie groups and
symmetric spaces.
Graduate Studies in Mathematics \bf 34\rm, AMS (2001)

\bibitem{K} Koosis, P.: 
Introduction to $H^p$-spaces.
Cambridge Tracts in Mathematics \bf 115\rm, Cambridge University
Press 1998.

\bibitem{MM} Matsushima, Y., Morimoto, A.:
Sur certaines espaces fibr\'es
holomorphes sur une vari\'et\'e de Stein.
\sl Bull.~Soc.~Math. France \bf 88\rm, 137-155 (1960)

\bibitem{R} Rudin, W.: 
Function theory in the unit ball of $\C^n$.
Grundlehren der math. Wiss. \bf 241\rm.
Springer 1980.

\bibitem{SZ} Saerens, R.; Zame, W.R.:
The Isometry Groups of Manifolds and the Automorphism Groups of Domains.
\sl Trans. A.M.S. \bf 301\rm, no. 1, 413-429 (1987)

\bibitem{St} Stein, K.: \"Uberlagerungen holomorph--vollst\"andiger komplexer
R\"aume. \sl Arch. Math. \bf VIII \rm, 354--361
(1956)

\bibitem{TS} Tumanov, A.E.; Shabat, G.B.:
Realization of linear Lie Groups by Biholomorphic Automorphisms
\sl Funct. Anal. Appl. \bf 24\rm, 255--257 (1991)

\bibitem{W1} Winkelmann, J.:
Invariant Hyperbolic Stein Domains.
\sl manu. math. \bf 79\rm, 329--334 (1993)

\bibitem{W2} Winkelmann, J.:
Realizing Countable Groups as Automorphism Groups of
Riemann Surfaces.
\sl documenta math.\rm, (2002)

\end{thebibliography}
\end{document}